\documentclass[11pt]{amsart}
\usepackage{graphicx}
\usepackage{epsfig}
\usepackage{amsmath}
\usepackage{stackrel, amssymb}
\usepackage{amscd}
\usepackage{latexsym}
\usepackage{tabularx}
\usepackage{a4wide}
\usepackage[usenames]{color}
\usepackage{subfigure}
\usepackage{url}
\usepackage{verbatim}
\usepackage{mathtools}
\usepackage{circuitikz}
\usepackage{mathabx}
\usepackage{tikz-cd}
\usepackage[T1]{fontenc}
\DeclareMathAlphabet{\mathdutchcal}{U}{dutchcal}{m}{n}

\newcommand{\Gm}{\ssub{\ssub[-1pt]{\G}!_{\mathbf m}}!}

\usepackage{enumitem} 
\setlist[enumerate]{label=\textnormal{(\arabic*)}}

\ctikzset{bipoles/resistor/height=0.15}
\ctikzset{bipoles/resistor/width=0.4}
\usepackage[
colorlinks, citecolor=blue,
pdfauthor={Omid Amini, Eduardo Esteves, Eduardo Garcez}, 
pdftitle={the-residue-polytope},
pdfstartview ={FitV},
]{hyperref}
\usepackage{tikz, float} 
\usetikzlibrary {positioning}
\usetikzlibrary{calc,decorations.markings}
\usetikzlibrary{shapes,snakes}
\numberwithin{equation}{section}
\tikzstyle{Cwhite}=[scale = .8,circle, fill = white, minimum size=3mm] 
\tikzstyle{Cgray}=[scale = .4,circle, fill = gray, minimum size=3mm] 
\tikzstyle{Cblack2}=[scale = .4,circle, fill = black, minimum size=5mm] 
\tikzstyle{Cblack}=[scale = .7,circle, fill = black, minimum size=3mm]
\tikzstyle{C0}=[scale = .9,circle, fill = black!0, inner sep = 0pt, minimum size=3mm]
\tikzstyle{C1}=[scale = .7,circle, fill = black!0, inner sep = 0pt, minimum size=3mm]
\tikzstyle{Cred}=[scale = .4,circle, fill = red, minimum size=3mm] 
\usepackage[matrix,arrow,tips,curve]{xy}
\usepackage{pb-diagram,pb-xy}
\usepackage{verbatim}
\usepackage{mathrsfs}
\usepackage{color}
\usepackage[leqno]{amsmath}
\usepackage{scalerel}
\usepackage{marginnote}

\newtheorem{thm}{Theorem}[section]

\newtheorem{lemma}[thm]{Lemma}

\newtheorem{prop}[thm]{Proposition}

\theoremstyle{definition}

\newtheorem{remm}[thm]{Remark}
\newenvironment{remark}
  {\pushQED{\qed}\remm}
  {\popQED\endremm}

\numberwithin{equation}{section}

\newcommand{\G}{\mathbf G}

\renewcommand{\k}{\ssub{\mathbf k}!}
\newcommand{\R}{\ssub{\mathbb R}!}
\renewcommand{\P}{\ssub[-2pt]{\mathbf P}!}
\newcommand{\RP}{\ssub[-2pt]{\mathbf P}!_{\res}}
\newcommand{\Q}{\ssub{\mathbf Q}!}

\newcommand{\Z}{\ssub{\mathbb Z}!}

\newcommand{\C}{{\mathscr C}}

\newcommand{\g}{\mathfrak{g}}
\newcommand{\red}{\ssub{\mathrm{red}}!}
\newcommand{\res}{\ssub{\mathrm{res}}!}
\newcommand{\Res}{\ssub{\mathrm{Res}}!}

\usepackage{todonotes}

\NewDocumentCommand{\ssub}{O{0pt} O{.9} m t! e{_^}}{
  #3%
  \IfValueT{#5}{
    \IfBooleanTF{#4}{\sb{\hspace{#1}\scaleobj{#2}{#5}}}{\sb{#5}}
  }
  \IfValueT{#6}{
  \IfBooleanTF{#4}{\sp{\hspace{#1}\scaleobj{#2}{#6}}}{\sp{#6}}
}
}

\NewDocumentCommand{\ssubb}{O{0pt} O{0pt} O{.9} m t! e{_^}}{
  #4%
  \IfValueT{#6}{
    \IfBooleanTF{#5}{\sb{\hspace{#1}\scaleobj{#3}{#6}}}{\sb{#6}}
  }
  \IfValueT{#7}{
  \IfBooleanTF{#5}{\sp{\hspace{#2}\scaleobj{#3}{#7}}}{\sp{#7}}
}
}
\ExplSyntaxOn
\NewDocumentCommand{\tossub}{o o m}{
  \expandafter\let\csname old\cs_to_str:N #3\endcsname#3
  \renewcommand#3%
  {\ssub[#1][#2]{\csname old\cs_to_str:N #3\endcsname}}
}
\ExplSyntaxOff

\newcommand{\stable}{\mathscr A}

\newcommand{\rquot}[2]{#1\bigl/#2}

\newcommand{\cd}{\mathit{cd}}
\newcommand{\dd}{\mathit{d}}

\newcommand{\VS}{\ssub{\mathbf U}!}

\tossub\omega
\tossub\Omega
\tossub\alpha
\tossub\beta
\tossub\sigma
\tossub\mu
\tossub\nu
\tossub\varphi
\tossub\delta
\tossub\Delta

\tossub\Gamma
\tossub\gamma
\tossub\zeta
\tossub\eta
\tossub\stable
\tossub\C
\tossub\iota
\tossub\rho
\tossub\varrho

\tossub\g
\tossub\red
\tossub\cd
\tossub\dd
\tossub\varpi
\tossub\lambda
\tossub\epsilon
\tossub\varepsilon
\tossub\theta
\tossub\xi

\newcommand{\ssM}{\ssub{M}!}

\newcommand{\ssI}{\ssub{I}!}

\newcommand{\rmH}{\ssub{\mathrm{H}}!}

\newcommand{\sspi}{\ssub{\pi}!}

\newcommand{\ssW}{\ssub[-1pt]{\mathrm W}!}
\newcommand{\ssn}{\ssub{n}!}
\newcommand{\ssq}{\ssub{q}!}
\newcommand{\ssh}{\ssub{h}!}

\newcommand{\ssE}{\ssub{E}!}

\newcommand{\ssd}{\ssub{d}!}
\newcommand{\ssx}{\ssub{x}!}
\newcommand{\ssy}{\ssub{y}!}
\newcommand{\ssz}{\ssub{z}!}

\newcommand{\ssc}{\ssub{\scaleto{\mathfrak{c}}{6pt}}!}

\newcommand{\ssm}{\ssub{m}!}
\newcommand{\ssp}{\ssub{p}!}
\newcommand{\ssF}{\ssub{F}!}
\newcommand{\ssS}{\ssub{S}!}

\newcommand{\filt}{\ssub{\mathrm{F}}!}

\renewcommand{\setminus}{\smallsetminus}
\renewcommand{\emptyset}{\varnothing}

\newcommand{\E}{\ssub{\mathbb E}!}
\newcommand{\A}{\ssub{A}!}   
\newcommand{\barA}{\ssub{\bar A}!}   

\newcommand{\he}{\ssub{\scaleto{\mathfrak{h}}{8pt}}!}  
\newcommand{\te}{\ssub{\scaleto{\mathfrak{t}}{6pt}}!}  

\makeatletter
\newsavebox\myboxA
\newsavebox\myboxB
\newlength\mylenA

\newcommand*\overbar[2][0.75]{%
    \sbox{\myboxA}{$\m@th#2$}%
    \setbox\myboxB\null
    \ht\myboxB=\ht\myboxA%
    \dp\myboxB=\dp\myboxA%
    \wd\myboxB=#1\wd\myboxA
    \sbox\myboxB{$\m@th\overline{\copy\myboxB}$}
    \setlength\mylenA{\the\wd\myboxA}
    \addtolength\mylenA{-\the\wd\myboxB}%
    \ifdim\wd\myboxB<\wd\myboxA%
       \rlap{\hskip 1\mylenA\usebox\myboxB}{\usebox\myboxA}%
    \else
        \hskip -0.5\mylenA\rlap{\usebox\myboxA}{\hskip 0.5\mylenA\usebox\myboxB}%
    \fi}

\newcommand{\comp}[1]{\overbar[.5]{#1}} 
\newcommand{\mg}{\mathscr M}
\newcommand{\mgbar}{\ssub[-1pt]{\comp\mg}!}

\newcommand{\proj}[1]{\theta_{_{ \hspace{-.07cm}#1}}}

\newcommand{\st}{\bigm|} 

\newcommand{\Fun}{\mathscr E} 
\tossub\Fun 

\newcommand{\id}{\mathrm{Id}}
\tossub\id
\newcommand{\zero}{0}
\newcommand{\one}{\ssub{\mathbf 1}!}
\tossub\zero
\tossub\one

\newcommand{\sspsi}{\ssub{\psi}!}

\tossub \chi

\newcommand{\varC}{\mathbf C}

\newcommand{\varR}{\mathbf R}

\newcommand{\ssA}{\ssub{A}!}
\newcommand{\ssV}{\ssub{V}!}
\newcommand{\ssH}{\ssub{H}!}

\newcommand{\summit}{\ssub{{\scaleto{\mathrm{s}}{5.5pt}}}!}
\newcommand{\ired}{\mathrm{ir}}

\newcommand{\codim}{\mathrm{codim}}
\newcommand{\ssu}{\ssub{u}!}

\definecolor{cadmiumgreen}{rgb}{0.0, 0.42, 0.24}
\definecolor{darkred}{rgb}{.6,0,0}
\definecolor{byzant}{rgb}{0.74, 0.2, 0.64}
 \definecolor{pblue}{rgb}{0.11, 0.22, 0.73}
\definecolor{pgreen}{rgb}{0.0, 0.65, 0.58}
\definecolor{aqua}{rgb}{0.0, 1.0, 1.0}
\definecolor{lblue}{rgb}{0.0, 0.55, 1.0}
\definecolor{pblue}{rgb}{0.11, 0.22, 0.73}

\DeclareFontFamily{U}{MnSymbolA}{}
\DeclareSymbolFont{MnSyA}{U}{MnSymbolA}{m}{n}
\DeclareFontShape{U}{MnSymbolA}{m}{n}{
    <-6>  MnSymbolA5
   <6-7>  MnSymbolA6
   <7-8>  MnSymbolA7
   <8-9>  MnSymbolA8
   <9-10> MnSymbolA9
  <10-12> MnSymbolA10
  <12->   MnSymbolA12}{}
\DeclareMathSymbol{\dashedleftarrow}{\mathrel}{MnSyA}{98}
\DeclareMathSymbol{\dashedrightarrow}{\mathrel}{MnSyA}{96}

\makeatletter
\newcommand{\topreleft}[1]{%
  \vbox {\m@th\ialign{##\crcr
  \topreleftfill \crcr
  \noalign{\kern-\p@\nointerlineskip}
  $\hfil\displaystyle{#1}\hfil$\crcr}}}

\newcommand{\topreright}[1]{%
  \vbox {\m@th\ialign{##\crcr
  \toprerightfill \crcr
  \noalign{\kern-\p@\nointerlineskip}
  $\hfil\displaystyle{#1}\hfil$\crcr}}}

\def\topreleftfill{%
  $\m@th%
  \dashedleftarrowtip%
  \mkern-1mu%
  \xleaders\hbox{$\mkern2mu\shortbar\mkern-1mu$}\hfill%
  \mkern1mu%
  \shortbar%
  \mkern0.5mu%
$}

\def\toprerightfill{%
  $\m@th%
  \mkern.5mu%
  \shortbar%
  \mkern-1mu%
  \xleaders\hbox{$\mkern2mu\shortbar\mkern-1mu$}\hfill%
  \mkern1mu%
  ‌​\dashedrightarrowtip%
$}

\def\dashedleftarrowtip{%
  \raisebox{\z@}[4.0pt][0.0pt]{$\mathord{\dashedleftarrow}$}}

\def\dashedrightarrowtip{%
  \raisebox{\z@}[4.0pt][0.0pt]{$\mathord{\dashedrightarrow}$}}

\def\shortbar{%
  \smash{\scalebox{0.4}[1.0]{$-$}}}
\makeatother

\newcommand{\ssvarphi}{\ssub{\varphi}!}

\newcommand{\bV}{\ssubb[-2pt]{\mathbf{\Upsilon}}!}

\newcommand{\LRC}{\ssub{{\mathcal T}_{\scaleto{\mathrm{loc}}{4pt}}}!}  

\newcommand{\Ros}{\ssub{{\mathcal T}_{\scaleto{\mathrm{Ros}}{4pt}}}!}  
\newcommand{\RosSp}{\ssub{\mathbfcal R}!}

\newcommand{\Glob}{\ssub{{\mathcal T}_{\scaleto{\mathrm{glob}}{5pt}}}!}

\DeclareMathAlphabet\mathbfcal{OMS}{cmsy}{b}{n}

\newcommand{\GlobSp}{\ssub{\mathbfcal{G}}!}

\newcommand{\mC}{\ssub{\mathrm{C}}!}
\newcommand{\rmC}{\mathrm{C}}
\newcommand{\lvl}{\ssub{L}!}

\newcommand{\levcom}{\ssub{\mathcal L\mathcal C}!}

\newcommand{\concom}{\ssub{\mathcal C\mathcal C}!}

\newcommand{\ssK}{\ssub{K}!}

\newcommand{\grass}{\mathrm{Grass}}

\newcommand{\subface}{\ssub{\prec}!}

\newcommand{\supface}{\ssub[-1pt]{\succ}!}
\newcommand{\supfaceq}{\ssub{\succeq}!}

\title{Residue polytopes}
 \author{Omid Amini}
 \address{CNRS - CMLS, \'Ecole Polytechnique, Palaiseau, France}
\email{omid.amini@polytechnique.edu}

 \author{Eduardo Esteves}
 \address{Instituto Nacional de Matem\'atica Pura e Aplicada, Estrada Dona Castorina 110,
22460-320 Rio de Janeiro RJ, Brazil}
\email{esteves@impa.br}

 \author{Eduardo Garcez}
 \address{Universidade Estadual do Cear\'a, Aracati, Brazil}
\email{jemgarcez@gmail.com}
 \date{\today}

\begin{document}
\begin{abstract} 
 A level graph is the data of a pair $(G,\pi)$ consisting of a finite graph $G$ and an ordered partition $\pi$ on the set of vertices of $G$. To each level graph on $n$ vertices we associate a polytope in $\R^n$ called its \emph{residue polytope}. We show that residue polytopes are compatible with each other in the sense that if $\pi'$ is a coarsening of $\pi$, then the polytope associated to $(G,\pi)$ is a face of the one associated to $(G,\pi')$. Moreover, they form all the faces of the \emph{residue polytope of $G$}, defined as the polytope associated to the level graph $(G, \sspi_0)$ for the trivial ordered partition $\sspi_0$.  The results are used in a companion work to describe limits of spaces of Abelian differentials on families of Riemann surfaces approaching a stable Riemann surface on the boundary of the moduli space.
\end{abstract}
\maketitle

\setcounter{tocdepth}{1}
\tableofcontents

\section{Introduction}
 Let $G=(V,E)$ be a fixed finite graph with set of vertices $V$, set of edges of $E$, and $\ssc$ connected components. For each positive integer $r$, let $[r]\coloneqq\{1,\dots, r\}$.

 A \emph{partition} of $V$ is a collection of nonempty pairwise disjoint subsets of $V$ whose union is $V$. A \emph{level structure} on $G$ is the data $\pi = (\sspi_1, \dots, \sspi_r)$ of an ordered partition of the vertex set $V$, that is, a partition of $V$ endowed with a total order on its set of elements.  It gives rise to a surjection $h=\ssh_\pi\colon V \to [r]$, taking each vertex to the part it belongs to. We refer to either $(G, \pi)$ or $(G,h)$ as a \emph{level graph}, and call $h$ the associated \emph{level function}. For $v\in V$, $h(v)$ is called the \emph{level} of $v$. For each pair of vertices $u,v\in V$, we write $u\subface_\pi v$ if $h(u) < h(v)$. We also denote by $\Pi=\Pi(G)$ the collection of ordered partitions of the vertex set $V$.

 We denote by $\E$ the set of \emph{arrows} of $G$: each edge $e=\{v,w\}$ gives rise to two arrows $vw$ and $wv$  in $\E$, written graphically, $v\to w$ and $w\to v$, respectively. If~$v=w$, we still have two arrows $v\leftrightarrow w$, corresponding to the two half edges issued from $e$. For the arrow $a=v \to w$, we call  $v$ the \emph{tail} and $w$ the \emph{head} of $a$, and write $\te_a=v$ and $\he_a=w$. The arrow in $\E$ with the reverse orientation is denoted $\bar a$; it has tail $w$ and head $v$.

A \emph{vertical edge} in a level graph $(G,\pi)$ is an edge $e=\{u,v\}$ with $\ssh_\pi(u)\neq \ssh_\pi(v)$. We denote the set of vertical edges by $\ssE_\pi$. A nonvertical edge is called \emph{horizontal}. They form the elements of the complement $\ssE_\pi^c = E\setminus \ssE_\pi$.  We denote by $\E_\pi$ the subset of $\E$ consisting of arrows on vertical edges. An element of $\E_\pi$ is called a vertical arrow. The complement $\E_\pi^c = \E \setminus \E_\pi$ is the set of horizontal arrows.

\begin{figure}[ht]
    \centering
\begin{tikzpicture}

\draw[line width=0.4mm] (5,-2) -- (6,0);
\draw[line width=0.4mm] (7,-2) -- (6,0);
\draw[line width=0.4mm] (7,-2) -- (7,0);
\draw[line width=0.4mm] (5,-2) -- (7,0);
\draw[line width=0.4mm] (8.5,-2) -- (7,0);
\draw[line width=0.4mm] (7,-2) -- (8.5,-2);
\draw[line width=0.4mm] (8.5,-2) -- (7,-2);

\draw[line width=0.1mm] (6,0) circle (0.8mm);
\filldraw[aqua] (6,0) circle (0.7mm);
\draw[line width=0.1mm] (7,0) circle (0.8mm);
\filldraw[aqua] (7,0) circle (0.7mm);
\draw[line width=0.1mm] (5,-2) circle (0.8mm);
\filldraw[aqua] (5,-2) circle (0.7mm);
\draw[line width=0.1mm] (7,-2) circle (0.8mm);
\filldraw[aqua] (7,-2) circle (0.7mm);
\draw[line width=0.1mm] (8.5,-2) circle (0.8mm);
\filldraw[aqua] (8.5,-2) circle (0.7mm);

\draw[black] (6,0.25) node{$\ssu_4$};
\draw[black] (7,0.25) node{$\ssu_5$};
\draw[black] (5,-2.3) node{$\ssu_1$};
\draw[black] (7,-2.3) node{$\ssu_2$};
\draw[black] (8.5,-2.3) node{$\ssu_3$};

\end{tikzpicture} 
\caption{A level graph $(G,\pi)$ with two levels, $\pi=(\sspi_1, \sspi_2)$ with $\sspi_1=\{\ssu_4, \ssu_5\}, \sspi_2 = \{\ssu_1, \ssu_2, \ssu_3\}$. The arrows $\ssu_1\ssu_4, \ssu_1\ssu_5, \ssu_2\ssu_4, \ssu_2\ssu_5, \ssu_3\ssu_5$ are upward. The edge $\{\ssu_2,\ssu_3\}$ is horizontal.}
\label{fig:digraph}
\end{figure}
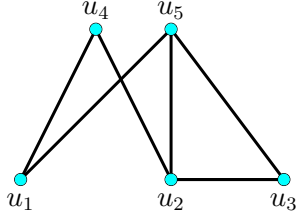

An arrow $u\rightarrow v$ in $\E_\pi$ is said to be compatible with the level structure if $u \supface_\pi v$. (Visually, both $\rightarrow$ and $\supface_\pi$ point to $v$.)   We denote by $\A_\pi$ the set of these arrows and refer to them as \emph{upward arrows}. Its complement in $\E_\pi$ is denoted $\barA_\pi$, and its elements are called \emph{downward arrows}. We have a bipartition $\E_\pi = \A_\pi \sqcup \barA_\pi$. Figure~\ref{fig:digraph} shows a level graph $G$ with two levels. In figures, we draw horizontal edges horizontally and vertical edges vertically. As the directions of the upward and downward arrows are visually clear, we drop them.

Let $\k$ be a field and set $\bV\coloneqq\k^{\E}$. For each $\psi \in \bV$ and $a\in \E$, denote by $\psi_a$ the $a$-coordinate of $\psi$. Notice that we can decompose $\bV =\bigoplus_{v\in V} \bV_v$ with $\bV_v\coloneqq\k^{\E_v}$  for each $v\in V$, with $\E_v$ the set of arrows $vw$ in $\E$ with tail $v$.

\smallskip

Given a level graph $(G,\pi =(\sspi_1, \dots, \sspi_r))$ with level function $h$, we define the \emph{residue space of $(G,\pi)$}, or simply \emph{$\pi$-residue space} if $G$ is fixed, as the subspace of $\bV$ denoted by $\GlobSp_\pi$  consisting of those elements $\psi$ which verify \emph{the residue conditions}~\ref{cond:R1} through \ref{cond:R4} listed as follows (see Section~\ref{sec:geometric-background} below for a discussion of their geometric origin):

\begin{enumerate}[label=(R\arabic*)]
    \item \emph{Vanishing along downward arrows}\label{cond:R1}
    \[\sspsi_a=0  \qquad \text{for every arrow $a \in\barA_{\pi}$}.\] 
    
    \smallskip
    
    \item \emph{Local residue conditions}\label{cond:R2}
    \[\sum_{a\in\E_v}\sspsi_a=0 \qquad \text{for every $v\in V$}.\]

    \smallskip
    \item \emph{Rosenlicht conditions}\label{cond:R3}
    \[
    \sspsi_a+\sspsi_{\bar a}=0 \qquad \text{ for each horizontal arrow $a$ for $\pi$}.  
    \]
    
    \smallskip
    \item \emph{Global residue conditions}\label{cond:R4}
    \[
    \sum_{a\in \ssA^{\Xi}_{\pi, n}}\psi_a=0
    \]
\end{enumerate}
for each level $n\in [r]$ and each connected component $\Xi$ of the subgraph $G[\ssV_{h< n}]$ of $G$ induced on the set of vertices $\ssV_{h< n} \subseteq V$ of level smaller than $n$, with $\A_{\pi,n}^{\Xi}\subseteq\A_{\pi}$ denoting the set of upward arrows with tail of level $n$ and head in $\Xi$. 

In Figure~\ref{fig:digraph}, taking $n=2$, the set $\ssV_{h< 2}$ consists of two vertices, $\ssu_4$ and $\ssu_5$, with the induced graph $G[\ssV_{h< 2}]$ consisting of no edges. Taking $\Xi$ to be the connected component $G[\ssV_{h< 2}]$ with vertex set $\{\ssu_5\}$, the set $\ssA^\Xi_{\pi,2}$ has three arrows, $\ssu_1\ssu_5, \ssu_2\ssu_5, \ssu_3\ssu_5$. The corresponding global residue condition is the equation $\psi_{\ssu_1\ssu_5}+\psi_{\ssu_2\ssu_5}+\psi_{\ssu_3\ssu_5}=0$.

Note that if $\sspi_0 =\{V\}$ is the trivial partition of $V$, Condition~\ref{cond:R4} is vacuum. The $\sspi_0$-residue space $\GlobSp_{\sspi_0}$ can be identified with the first homology group $\ssH_1(G, \k)$ (equivalently, the space of $\k$-valued flows in the graph), which has dimension equal to the \emph{genus} $g=g(G)$ of $G$ given by $g\coloneqq |E|-|V|+\ssc$. Our first theorem states that this holds for all residue spaces.

\begin{thm}\label{thm:dimG} Let $\pi$ be an ordered partition of $V$. We have 
\[
\dim_{\k}\GlobSp_\pi=|E|-|V|+\ssc.
\]
\end{thm}

To each residue space $\GlobSp_{\pi}$ for $\pi\in\Pi$, we associate a polytope $\P_{\pi} \subseteq \R^V$ as follows. For each subset $I\subseteq V$, consider the projection map 
\[
\proj{I} \colon \bigoplus_{v\in V} \bV_v \to \bigoplus_{v\in I} \bV_v.
\]
The $\pi$-residue space $\GlobSp_{\pi}$ lives in $\bV= \bigoplus_{v\in V} \bV_v$. Define the function  
\[
\gamma!_{\pi}\colon\ssub{2}!^V \to \Z, \qquad\qquad \dim_{\k} \gamma!_{\pi}(I) \coloneqq \proj{I}(\GlobSp_{\pi}) \quad \text{for each $I\subseteq V$}.
\]
This is a \emph{submodular function}, meaning that it verifies $\gamma!_\pi(\emptyset)=0$ and the following set of inequalities
\[
\gamma!_{\pi}(I) + \gamma!_{\pi}(J) \geq \gamma!_{\pi}(I\cup J) +\gamma!_{\pi}(I\cap J) \qquad \text{for all } I, J\subseteq V.
\]
We then define $\P_\pi$ as the base polytope of the polymatroid defined by $\gamma!_{\pi}$, namely,
\[
\P_\pi \coloneqq \Bigl\{q\in \R^V \, \st\, q(I) \leq \gamma!_\pi(I) \text{ for all } I\subseteq V \Bigr\}.
\]
Here, for $q\in \R^V$ and $I\subseteq V$,  we set $q(I) = \sum_{v\in I}q(v)$. Then, we have the inclusion $\P_\pi \subseteq \Delta!_g$ where $\Delta!_g \subseteq \R_{\geq 0}^V$ is the standard simplex of size $g$, consisting of those points $q\in \R^V$ with $q(v)\geq 0$ for all $v\in V$ and $q(V)=g$.

For the trivial ordered partition $\sspi_0=\{V\}$, we have 
\[
\gamma!_{\sspi_0}(I) = g - g(\ssI^c) \qquad \textrm{ for all } I \subseteq V
\]
where $\ssI^c = V\setminus I$, and $g(\ssI^c)$ is the genus of the induced subgraph $G[\ssI^c]$. Equivalently, $\gamma!_{\sspi_0}(I)$ is the genus of the graph obtained from $G$ by contracting each connected component of $G[\ssI^c]$ into a vertex. We call $\P_{\sspi_0}$ the \emph{residue polytope of $G$} and denote it by $\RP=\RP(G)$.

\smallskip

Our second theorem is stated as follows. Given two ordered partition $\pi=(\sspi_1, \dots, \sspi_r)$ and $\pi'=(\sspi'_1,\dots, \sspi'_{s}) \in \Pi$, we say $\pi'$ is a \emph{coarsening} of $\pi$ and write $\pi' \supfaceq \pi$ if each element of $\pi$ is a subset of an element of $\pi'$, and the total orders are compatible. More precisely, we have $\pi' \supfaceq \pi$ if the surjection $h'\colon \ssV \to [s]$ associated to $\sspi'$ factors through the surjection $h\colon \ssV\to [r]$ associated to $\pi$, and the induced map  $c\colon[r]\to[s]$ satisfies $c(\ssn_1)\leq c(\ssn_2)$ for each $\ssn_1,\ssn_2\in[r]$ with $\ssn_1\leq \ssn_2$. 
In particular, note that $\sspi_0$ is a coarsening of any $\pi\in \Pi$. Our second main result is stated as follows.

\begin{figure}
    \centering
   \scalebox{.4}{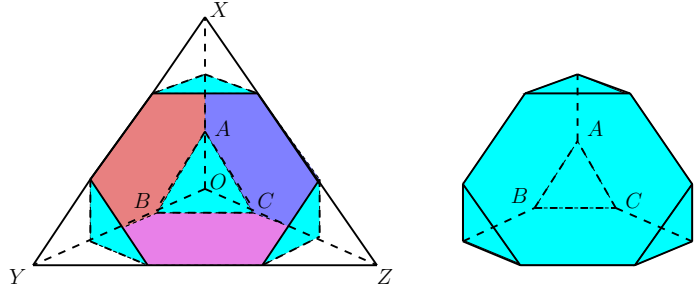}
    \caption{The residue polytope of the complete graph $K_4$ on four vertices on the right. The figure on the left shows the position of $\RP(K_4)$ within the simplex $OXYZ$ of width three in $\R^4$, with vertices $O=(3,0,0,0)$, $X=(0,3,0,0)$, $Y=(0,0,3,0)$ and $Z=(0,0,0,3)$, using the projection to $\R^3$ given by the last three coordinates.}
    \label{fig:residue-polytope-k4}
\end{figure}

\begin{thm}\label{thm:faces-residue-polytope} For each  $\pi\in\Pi$, the polytope $\P_{\pi}$ is a face of $\RP$. Furthermore, each face of $\RP$ is of this form. In addition, if $\pi' \supfaceq \pi$, then $\P_{\pi'} \supseteq \P_{\pi}$, and so $\P_{\pi}$ is a face of $\P_{\pi'}$.
\end{thm}
The residue polytope of the complete graph $\ssK_4$ on $4$ vertices in depicted in Figure~\ref{fig:residue-polytope-k4}.

In order to prove the first theorem above, we define for each $\pi\in \Pi$, a flag of subspaces
\[
\GlobSp_\pi \subseteq \RosSp_\pi \subseteq \bV^0_\pi \subseteq \bV_\pi \subseteq \bV
\]
by imposing Conditions~\ref{cond:R1} through~\ref{cond:R4}, one after the other. Namely, $\bV_\pi \subseteq \bV$ is obtained by imposing \ref{cond:R1},
\begin{equation}\label{eq:flag1}
\bV_\pi \coloneqq \bigl\{\psi\in \bV \, \st \, \psi_a =0 \textrm{ for every }  a \in \barA_{\pi} \bigr\}\subseteq\bV.
\end{equation}
Further down, the subspace $\bV^0_\pi \subseteq \bV_\pi$ is obtained by imposing \ref{cond:R2},
\begin{equation}\label{eq:flag2}
\bV_{\pi}^0 \coloneqq \bigl\{\psi \in \bV_\pi \,\st\, \sum_{a\in \E_v} \psi_a =0  \text{ for each $v\in V$}\bigr\}.
\end{equation}
Finally, the subspace $\RosSp_\pi \subseteq \bV_\pi^0$ is defined by imposing \ref{cond:R3},
\begin{equation}\label{eq:flag3}
\RosSp_\pi\coloneqq \bigl\{\psi\in\bV^0_{\pi}\,\st\,\psi_a+\psi_{\bar a}=0\quad \text{for all } a\in\E\setminus\E_{\pi}\bigr\}.
\end{equation}
We obtain the dimensions of the spaces appearing in the flags in Sections~\ref{sec:R1-2-3} and \ref{sec:grc}, which lead to the proof of Theorem~\ref{thm:dimG}.

\smallskip

In order to prove the second theorem, we show that for each pair of ordered partitions $\pi, \pi'\in\Pi$ with $\pi' \supfaceq \pi$, the $\pi$-residue space $\GlobSp_\pi$ can be obtained as a limit of a family of subspaces of $\bV$ with the same associated submodular function as $\gamma!_{\pi'}$. More precisely, let $\varR = \k[t,t^{-1}]$ be the ring of Laurent polynomials with $\k$-coefficients and denote by $\Gm^V(\varR)$ the set of invertible elements in $\varR$. An integer valued function on $V$ gives rise to an ordered partition of $V$ induced by its level sets. For each integer-valued function $d\colon V\to\Z$, $v\mapsto \ssd_v$, with induced ordered partition $\pi$, let $x\in\Gm^V(\varR)$ be given by $\ssx_v\coloneqq t^{-\ssd_v}$ for each $v\in V$, and put
\[
\ssub{(x\cdot\psi)}_{a}\coloneqq
\ssx_{v}\sspsi_a=t^{-\ssd_v}\sspsi_a \qquad \text{for each $\psi\in\bV$ and each $v\in V$ and $a\in\E_v$}.
\]
This defines an action of $\Gm^V(\varR)$ on $\bV\otimes_{\k} \varR$. By studying this action and using Theorem~\ref{thm:dimG}, we show in Theorem~\ref{thm:limGlobSp} that if $\pi'\supfaceq \pi$, then
\[
\lim_{t\to 0}x\cdot\GlobSp_{\pi'}=\GlobSp_{\pi},
\]
and deduce that $\GlobSp_{\pi}$ is the \emph{splitting} of $\GlobSp_{\pi'}$ with respect to $\pi$, in the sense of Section~\ref{sec:splitting} (see Proposition~\ref{prop:splitting} and~\eqref{eq:splitting2}).

Combining this with Proposition~\ref{prop:splitting}, we prove that the submodular function $\gamma!_{\pi}$ is a \emph{splitting} of the submodular function $\gamma!_{\pi'}$ in the sense of~\cite[\S2.4]{AE-modular-polytopes}, see Section~\ref{sec:upper-semicontinuity}. Using the characterization of faces of polymatroids, we finish the proof of Theorem~\ref{thm:faces-residue-polytope}.

Notice that our Theorem~\ref{thm:limGlobSp} provides a new perspective on global residue conditions, by showing that they arise as degenerations of the Rosenlicht and local residue conditions.

\subsection{Abelian differentials on Riemann surfaces and their limits}  \label{sec:geometric-background}

The results of this paper have their source of motivation in our companion work \cite{AEG24}, in which we study the limits of spaces of Abelian differentials on Riemann surfaces. We briefly discuss the setup.

For a smooth compact Riemann surface $S$ of genus $g$, the space of Abelian differentials on $S$, denoted by $\ssH_S$, is the complex vector space of dimension $g$ given by the global sections of the canonical sheaf $\omega!_S$. 
Given a stable Riemann surface $X$ of genus $g$ lying on the boundary of the Deligne-Mumford compactification $\mgbar_g$, we are interested in describing all the limits of the spaces $\ssH_{S_t}$ for any one-parameter family of smooth Riemann surfaces approaching $X$ in the moduli space. This has been studied previously for very special $X$, compact type~\cite{EH87} and two-component~\cite{EM}, in the case the nodes $\ssp^a$ for $a\in \E_u$ are general on $\varC_u$ for each vertex $u\in V$.

The stable Riemann surface $X$ has a dual graph $G=(V,E)$ with each vertex $v\in V$ corresponding to a component $\varC_v$ in the desingularization of $X$, and each edge $e\in E$ corresponding to a node $\ssp^e$ on $X$. If $e$ connects $u$ and $v$, then $\ssp^e$ lies on $\varC_u$ and $\varC_v$.  For each arrow $a=u\to v$ in $\E$ over the edge $e=\{u,v\}$, denote by $\ssp^a$ the point on $\varC_u$ over the node $\ssp^e$ of $X$. 

The tropicalizaion of the family $S_t$ yields an edge length function $\ell\colon E \to \R_{>0}$. Moreover, as we showed in~\cite{AE-modular-polytopes}, tropicalization also associates a limit space $\ssW_h \subset \Omega\coloneqq \bigoplus_{v\in V}\Omega!_v$ for each  function $h\colon V \to \R$, where $\Omega!_v$ is the space of meromorphic differentials on  $C_v$, for each $v\in V$. In particular, we can view each element $\alpha \in \ssW_h$ as a collection of meromorphic differentials $\alpha!_v \in \Omega_v$ for $v\in V$. This way, we obtain a residue map 
\[
\Res \colon \ssW_h \to \bV = {\mathbb C}^{\E}, \qquad \qquad \alpha=(\alpha!_v)_{v\in V} \mapsto \Bigl[a=vu \mapsto \res_{\ssp^a}(\alpha!_v)\Bigr].
\]

Viewing each function $h$ as a level function gives rise to an ordered partition $\sspi_h$ of $V$. Let $\GlobSp_{\sspi_h}$ be the corresponding $\sspi_h$-residue space. Using results from \cite{BCGGM18, TT22} and~\cite{AE-modular-polytopes}, we prove in \cite[Thm 5.1]{AEG24} that the image of the residue map lies in $\GlobSp_{\sspi_h}$, that is, any element $\psi=\Res(\alpha)$ for $\alpha\in \ssW_h$ verifies Conditions~\ref{cond:R1} through \ref{cond:R4}. (The first three are classical, the global residue conditions were discovered in \cite{BCGGM18}.) 

The results of this paper play a crucial role in our description in \cite{AEG24} of all possible limits of $\ssH_{S_t}$, for all families $S_t$ as above approaching the stable Riemann surface $X$ with an arbitrary number of components, meeting in whatever ways, as long as the nodes on each component are in general position. Using them, we prove that there is a projective variety which parametrizes all these limits.

\section{Set-theoretically independent collections}\label{sec:stic} We start by formulating two basic linear-algebra results, frequently used in the sequel.

We fix a field $\k$. For a finite set $A$, we denote by $\k^A$ the $\k$-vector space of functions $\psi\colon A \to \k$, $a\mapsto \psi_a$.

 For each $\psi \in\k^A$, denote by $[\psi]$ the support of $\psi$ defined by 
 \[
 [\psi]\coloneqq \left \{a\in A\,\st \,\psi(a)\neq 0\right\}. 
 \] 
For a finite collection of vectors $\Psi\subset\k^A$, let $[\Psi]\subseteq A$ be the union of $[\psi]$ for $\psi\in\Psi$. We say that $\Psi$ is \emph{set-theoretically independent} if the $[\psi]$ for $\psi\in \Psi$ form a partition of $[\Psi]$. Obviously, if this happens, then $\Psi$ is linearly independent, and any subcollection $\Psi' \subseteq \Psi$ will be  set-theoretically independent as well.

\smallskip

We say that two set-theoretically independent collections of vectors 
$\Phi, \Psi \subset \k^A$ are \textit{related} when there are two nonempty subcollections  $\Phi' \subseteq \Phi$  and $\Psi' \subseteq \Psi$ that yield partitions of the same subset of $A$, that is, when $[\Phi'] = [\Psi']$. Otherwise, we call the two collections \emph{unrelated}. Note that in the latter case, $\Phi\cap \Psi=\emptyset$.

\begin{prop}\label{prop:unlinked}
    Let $\Phi, \Psi \subset\k^A$ be two unrelated set-theoretically independent collections of vectors. Then, the union $\Phi  \cup \Psi $ is linearly independent.
\end{prop}

\begin{proof} For the sake of contradiction, assume there exists a nontrivial relation
\begin{equation}\label{abli}
\sum_{\varphi \in \Phi} c_{\varphi}\varphi + \sum_{\psi \in \Psi}d_{\psi}\psi = 0.
\end{equation}
Since $\Phi$ and $\Psi$  are unrelated, we have 
\[
\bigcup_{\substack{\varphi:\, c_{\varphi}\neq0}} [\varphi] \neq  \bigcup_{\substack{\psi:\, d_\psi\neq0}} [\psi].
\]
We may assume w.l.o.g. that there is an $a\in A$ included in the left-hand side but not in the right-hand side. Then, there exists $\varphi\in\Phi$ with $\varphi(a)\neq 0$ and $c_{\varphi}\neq0$, unique with this property because $\Phi$ is set-theoretically independent. Thus, evaluating \eqref{abli} at $a$ gives $c_{\varphi}\varphi(a)=0$, a contradiction. 
\end{proof}

Two set-theoretically independent subsets 
$\Phi  ,\Psi \subseteq \k^A$ are called \emph{properly unrelated} if for each nonempty subcollections $\Phi'\subseteq \Phi$ and $\Psi'\subseteq \Psi$ with either $\Phi'\neq \Phi$ or $\Psi'\neq \Psi$, we have $[\Phi'] \neq [\Psi']$. 

\begin{prop}\label{prop:properlyunlinked} Let $\Phi,\Psi\subset\k^A$ be properly unrelated set-theoretically independent collections of vectors. Then, for each $\varphi\in\Phi$, the collection $ \left(\Phi \cup \Psi\right) \setminus \{\varphi\}$ is linearly independent.  
\end{prop}

\begin{proof} This follows from Proposition~\ref{prop:unlinked} applied to $\Phi\setminus \{\varphi\}$ and $\Psi$.
\end{proof}

\begin{remark} The above results generalize in an obvious way to the setup of a finite-dimensional vector space equipped with a basis.
\end{remark}

\section{Imposing the first three conditions} \label{sec:R1-2-3}
Let $\E$ be as before the set of arrows on edges of $G$, each edge $e=\{u,v\}$ giving rise to two arrows $uv$ and $vu$ in $\E$. For a subset $A \subseteq \E$, we denote by $\bar A$ the set consisting of the arrows $\bar a$ for all $a \in A$. Given subsets $I,J \subseteq V$, we denote by $\ssA_I$ the set of arrows in $A$ with tail in $I$, and denote by $A(I,J)$ the set of arrows in $A$ with tail in $I$ and head in $J$. If $I=J$, we simply write $A(I)$. Similarly, for $F\subseteq E$ and $I, J\subseteq V$, we define $\ssF_I$ as the set of edges of $F$ with one vertex in $I$, $\ssF(I,J)$ being the set of edges of $F$ with one vertex in $I$ and the other in $J$, and $\ssF(I)$ those edges of $F$ with both vertices in $I$.

We fix in this section and the next an ordered partition $\pi=(\sspi_1, \dots, \sspi_r)$ on $V$, and consider the level graph $(G, \pi)$. Let $h\colon V \to [r]$ be the corresponding level function. Let $\ssE_\pi$, $\ssA_\pi$, and $\barA_\pi$ be, respectively, the set of vertical edges,  upward arrows, and downward arrows in $(G,\pi)$. 

Given $n \in [r]$, a connected component of the graph $G[\sspi_n]$, induced on vertices of level $n$, is called a \emph{component of level $n$} or  a \emph{level $n$ component} of  $(G,\pi)$. Denote by $\levcom_n$ the collection of these level components. 

A \emph{summit} of $(G, \pi)$ is a level component $S$ containing no vertex which is the tail of an upward arrow. If $S$ is a singleton (one vertex, no edges), then the summit is called \emph{irreducible}; otherwise, it is called \emph{reducible}. We denote by $\summit_{\ired}$, resp.~$\summit_{\red}$, the number of irreducible, resp.~reducible, summits in $G$, and set  $\summit\coloneqq\summit_{\ired}+\summit_{\red}$, the total number of summits. 
For example, in Figure \ref{fig2}, the level components consisting of single vertices $\ssu_{5}, u_b$ and $u_c$ are the irreducible summits, whereas, the level components $G[\{\ssu_9 , u_a\}]$ and $G[\{\ssu_7 , \ssu_8\}]$ are the reducible summits. We thus have $\summit_{\ired}=3$ and $\summit_{\red}=2$; in total, there are $\summit=5$ summits.

We consider in the sequel the $\k$-vector space $\bV=\k^\E$ of $\k$-valued functions on $\E$.  Denote by $\langle\cdot,\cdot\rangle \colon \bV \times \bV \to \k$ the natural bilinear form.  For $v\in V$ and $S\subseteq V$, define
\[
\bV_v\coloneqq \k^{\E_v} \qquad \text{ and } \qquad \bV_S \coloneqq \bigoplus_{v\in S} \bV_v.
\]
For $a\in \E$, we denote by $\one_a$ the function on $\E$ that takes value 1 on $a$ and 0 elsewhere.  

\subsection{Vanishing residue along downward arrows} Let $\bV_\pi\subseteq\bV$ be the vector subspace given by the \emph{vanishing conditions along downward arrows}
\begin{equation}\label{eq:DAC}
\langle \one_a,\cdot\rangle=0 \qquad \qquad \text{for all } \,\, a\in\barA_\pi.
\end{equation}
Equivalently, $\bV_\pi$ is the set of $\psi\in\bV$ such that $\psi_a=0$ for each downward arrow $a\in\barA_\pi$. We have a decomposition $\bV_\pi=\bigoplus \bV_{\pi,v}$, where $\bV_{\pi,v}\coloneqq\bV_\pi\cap\bV_v$. The following is straightforward.
\begin{prop}\label{prop:DAC} We have $\dim\bV_\pi= 2|E| -|\ssE_\pi|$.
\end{prop}

For each $S\subseteq V$, let $\bV_{\pi,S}\coloneqq\bV_\pi\cap\bV_S$. Notice that $\bV_{\pi,S}$ has a natural basis consisting of the characteristic vectors $\one_a$ for all horizontal and upward arrows $a\in\E_S$. We will use this basis when applying the results in Section~\ref{sec:stic} to $\bV_{\pi,S}$.

\subsection{Local residue conditions} For each $v\in V$, let $\one_v \in \bV_\pi$ be the characteristic vector of $\E_v \setminus \barA_v$, that is, $\one_v \coloneqq  \sum_{a\in\E_v \setminus \barA} \one_a$. By definition, $\bV^0_\pi\subseteq\bV_\pi$ is the vector subspace defined by the \emph{local residue conditions}
\begin{equation}\label{eq:local-residue}
\langle \one_v,\cdot\rangle=0 \qquad \qquad \text{for all } \,\, v\in V.
\end{equation} 
Let $\bV^0_{\pi,v} \subseteq \bV_{\pi,v}$ be the subspace defined by
\[
\bV^0_{\pi,v}\coloneqq \Bigl\{\psi \in\bV_{\pi,v}\,\st\, 
\langle\one_v,\psi\rangle=0\Bigr\} =\Bigl\{\psi \in\bV_{\pi,v}\,\st\, 
\sum_{a\in\E_v\setminus \barA_{\pi,v}}\psi_a=\sum_{a\in\E_v}\psi_a=0\Bigr\}.
\]
We define $\bV^0_\pi\coloneqq \bigoplus_{v\in V}\bV^0_{\pi,v}$. Note that for each $v\in V$, $\bV_{\pi,v}^0=\bV_{\pi}^0\cap\bV_v$. For each $S\subseteq V$, we set $\bV_{\pi,S}^0\coloneqq\bV_{\pi}^0\cap\bV_S$.

\smallskip
Denote by  $\LRC\subseteq\bV_\pi$ the collection of all nonzero vectors $\one_v$, $v\in V$,
\[
\LRC \coloneqq \left\{\one_v \,\st\, v\in V,\, \one_v\neq 0\right\}.
\] 
Note that $\one_v=0$ if and only if $\E_v$ contains only downward arrows, that is, the singleton $v$ is an irreducible summit of the level graph $G$. Also, clearly, $\LRC$ is set-theoretically independent. Therefore, combined with Proposition~\ref{prop:DAC}, we get the following result.

\begin{prop}\label{prop:LRC} We have  $\codim(\bV^0_\pi,\bV_\pi)=|V|-\summit_{\ired}$.
\end{prop}

\subsection{Rosenlicht residue conditions}  For each horizontal edge $e \in \ssE_\pi^c=E\setminus \ssE_\pi$, let $\one_e \coloneqq \one_{a}+\one_{\bar a}$ be the characteristic vector of $e$, where $a$ and $\bar a$ are the two arrows on $e$.

Define the vector subspace $\RosSp_\pi \subseteq\bV^0_\pi$ by \emph{Rosenlicht residue conditions} 
\begin{equation}\label{eq:Rosenlicht}
\langle\one_e,\cdot\rangle=0 \qquad \qquad  \text{for all }\,\, e\in  \ssE_\pi^c.
\end{equation}
Let $\Ros\subseteq\bV_\pi$ be the collection of  characteristic vectors $\one_e$ for $e$ horizontal, that is,
\[
\Ros \coloneqq \bigl\{\one_e \,\st\, e\in \ssE_\pi^c \bigr\}.
\] 
This collection is set-theoretically independent in $\bV_\pi$. However, Equations~\ref{eq:Rosenlicht} are defined in the subspace $\bV_\pi^0$, and we need to do some extra work to get the dimension of $\RosSp_\pi$.

For each reducible summit with set of vertices $S \subseteq V$, each arrow in $\E_S$ is either horizontal or downward. Therefore, we have the equation
 \[
 \sum_{v\in S} \one_v = \sum_{e\in\ssE_{\pi}^c(S)} \one_e +\sum_{a\in\barA_{\pi,S}}\one_a,
 \]
 from which we deduce the following equation 
 \[
 \sum_{v\in S} \langle\one_v,\cdot\rangle =0 
 \]
on $\bV_\pi^0$. We will show below that these are the only relations the Rosenlicht conditions satisfy; see Lemma~\ref{lem:unrelated1}. From this, combined with Proposition~\ref{prop:LRC}, we will deduce the following. 

\begin{prop}\label{codimRos}
We have $\codim(\RosSp_\pi,\bV^0_\pi)=|\ssE_\pi^c|-\summit_{\red}$.
\end{prop}

In preparation for the proof, for each $n\in [r]$, consider the set $\sspi_n \subseteq V$ of vertices of level $n$, and the set $\levcom_n$  of level $n$ components. Since each element $\mC\in \levcom_n$ is entirely determined by its set of vertices, abusing  notation, we simply use $\mC$ for the set $V(\mC)$. We thus have 
 \[
 \sspi_n = \bigsqcup_{\rmC \in \levcom_n} \mC.
 \]
We can decompose 
\[
\bV^0_\pi=\bigoplus_{n\in [r]}\bV^0_{\pi,\sspi_n} \qquad \quad \textrm{and}\quad\qquad  \bV^0_{\pi,\sspi_n}=\bigoplus_{\rmC\in \levcom_n}\bV^0_{\pi,\rmC}.
\] 

Note that each vertex $v\in V$ belongs to a unique level component $\mC \in \levcom_n$ for unique $n\in [r]$, and we have $\one_v\in\bV_{\pi,\mC}$. Also, each horizontal edge $e\in \ssE_\pi^c$ connects a pair of vertices in a unique level component $\mC \in \levcom_n$, $n\in [r]$, in particular $\one_e\in\bV_{\pi,\mC}$. We thus get a decomposition 
\[
\RosSp_\pi = \bigoplus_{n\in [r], \, \mC\in \levcom_n} \RosSp_\pi \cap \bV_{\pi,\mC}.
\]
This implies that
\begin{equation}\label{eq:codim-Ros}
\codim(\RosSp_\pi,\bV^0_\pi)=\sum_{n\in [r],\,\rmC\in \levcom_n} \codim\Bigl(\RosSp_\pi \cap\bV_{\pi,\rmC}, \bV^0_{\pi,\rmC}\Bigr).
\end{equation}
For each $n \in [r]$ and $\mC\in \levcom_n$, define the two collections of vectors $\LRC^{\, n,\rmC}, \Ros^{\,n,\rmC} \subset \bV_{\pi,\rmC}$ by
\[
\LRC^{\, n,\rmC}\coloneqq \left\{\one_v \,\st v\in \rmC\right\}\qquad\text{and}\qquad
\Ros^{\,n,\rmC}\coloneqq\left\{\one_e\, \st\, e \in \ssE_\pi^c(\rmC)\right\}.
\]

\begin{lemma}\label{lem:unrelated1} Notation as above, for each $n\in [r]$ and $\mC\in \levcom_n$, the two collections $\LRC^{\,n,\rmC}$ and $\Ros^{\,n,\rmC}$ of vectors in $\bV_{\pi,\rmC}$ are properly unrelated. Moreover, they are related if and only if the level component $\mC$ is a reducible summit of the level graph $G$. 
\end{lemma}

Assuming this, we prove the proposition.

\begin{proof}[Proof of Proposition~\ref{codimRos}] We combine the above lemma with Propositions~\ref{prop:unlinked}~and~\ref{prop:properlyunlinked} applied to the two collections of vectors $\LRC^{\,n,\rmC}$ and $\Ros^{\,n,\rmC}$ in $\bV_{\pi,\rmC}$, for $n\in [r]$ and $ \mC\in \levcom_n$, to infer that 
\[
\codim\bigl(\RosSp_\pi \cap\bV_{\pi,\rmC},\bV^0_{\pi,\rmC}\big)=\left|\Ros^{\,n,\rmC}\right|-\epsilon(\rmC),
\]
where $\epsilon(\rmC)=1$ if $\mC$ is a reducible summit and $\epsilon(\rmC)=0$ otherwise. 
Combining Equation~\eqref{eq:codim-Ros} with the observation that
\[\sum_{n\in [r], \, \rmC\in \levcom_n} \left|\Ros^{\,n,\rmC}\right| = |\ssE_\pi^c|\]
 concludes the proof.
\end{proof}

\begin{proof}[Proof of Lemma~\ref{lem:unrelated1}] Let $\Phi \subseteq \LRC^{\,n,\rmC}$ and $\Psi \subseteq \Ros^{\,n,\rmC}$ be two subcollections of vectors with $[\Phi] = [\Psi]$, that is, 
\begin{equation}\label{eq:partition1}
   \bigcup_{\varphi \in \Phi}[\varphi] = \bigcup_{\psi \in \Psi}[\psi],
\end{equation} 
so that they yield partitions of the same subset of $\E \setminus \barA_\pi$.

We claim that $\Phi=\LRC^{\,n,\rmC}$. For the sake of  contradiction, suppose this is not the case and consider the proper subset $Z \subset \mC$ consisting of those vertices $u\in \rmC$ with  $\one_u\in \Phi$.  Since $\rmC$ is connected, there is an edge $e =\{u,v\}$ connecting a pair of vertices $u, v$ of $\rmC$ such that $u \in Z$ and $v \in \rmC \setminus Z$. This implies that $\one_{u}\in\Phi$ but $\one_v\not\in\Phi$. Note that $e$ is horizontal, i.e., $e\in\ssE_\pi^c$. Denote by $uv$ and $vu$ the two arrows on $e$. Since $uv\in [\one_{u}]$, the equality in \eqref{eq:partition1} implies that $uv\in [\psi]$ for some $\psi\in\Psi$. Necessarily, we have $\psi = \one_e = \one_{uv}+\one_{vu}$. But then, since $[\one_e] = \{uv, vu\}$, we also have $vu\in \bigcup_{\psi \in \Psi}[\psi]$, and hence $vu \in [\varphi]$ for some $\varphi\in\Phi$. This is possible only if $\varphi=\one_{v}$, a contradiction, proving the claim.

The union of the sets $[\varphi]$ for $\varphi\in\LRC^{\,n,\rmC}$ is $\E_{\mC}\setminus \barA_\pi$, which clearly contains $\E(\rmC)$, the set of (horizontal) arrows connecting vertices of $\rmC$. This set can be in the union of the $[\psi]$ for $\psi\in \Psi \subseteq \Ros^{\,n,\rmC}$ only if $\Psi=\Ros^{\,n,\rmC}$. We conclude that $\LRC^{\,n,\rmC}$ and $\Ros^{\,n,\rmC}$ are properly unrelated.

\smallskip

To prove the second statement, note that the two collections of vectors $\LRC^{\, n,\rmC}$ and $\Ros^{\,n,\rmC}$ are related if and only if they are both nonempty and we have 
 \[\bigcup_{\varphi \in \LRC^{\,n,\rmC}}[\varphi] = \bigcup_{\psi \in \Ros^{\,n,\rmC}}[\psi].\]
The above equality implies that all the arrows appearing in $[\varphi]$, $\varphi \in \LRC^{\,n,\rmC}$, are horizontal, that is, $\rmC$ is a summit. Moreover, if $\Ros^{\,n,\rmC}$ is nonempty, the summit must be reducible. This proves one direction of the statement. To finish, note that if $\rmC$ is a reducible summit, then $\LRC^{\,n,\rmC}$ and $\Ros^{\,n,\rmC}$ are nonempty, and $[\LRC^{\,n,\rmC}]=\E(\rmC) = [\Ros^{\,n,\rmC}]$.
\end{proof}

\section{Global residue conditions and proof of Theorem~\ref{thm:dimG}}\label{sec:grc} 
We keep the notation as in the previous section: $\pi=(\sspi_1, \dots, \sspi_r)$ is a level structure on $G$, and $h\colon V\to [r]$ is the corresponding level function. 

\smallskip

For each $n\in [r]$, let 
\[
\ssV_{h< n}\coloneqq \bigcup_{i< n} \sspi_i  \qquad \qquad \text{and}  \qquad \qquad \ssV_{h\leq n}\coloneqq \bigcup_{i\leq n} \sspi_i = \ssV_{h< n} \sqcup \sspi_n.
\] 
Let $\concom_{h< n}$ and $\concom_{h\leq n}$ be the set of connected components  of $G[\ssV_{h< n}]$ and  $G[\ssV_{h \leq n}]$, respectively.

\smallskip

A connected component $\Xi$ of $G[\ssV_{h< n}]$ is called \emph{special} if there is an arrow with tail in $\sspi_n$ and head in $\Xi$. Such an arrow is necessarily upward. Let  $\concom^*_{h< n} \subseteq \concom_{h< n}$ be the set of special connected components of $G[\ssV_{h< n}]$.  We also set $\concom^*_{h\leq n}\coloneqq \concom^*_{h< n+1}$. In Figure~\ref{fig2}, the connected component of the two upper levels formed by the vertices $\ssu_9$ and $u_a$ is special.   For each $\Xi \in \concom^*_{h<n}$,  let $\ssA_{\pi,n}^{\Xi}\subseteq \ssA_\pi$ be the set of upward arrows with tail in $\sspi_n$ and head in $\Xi$, namely,
\[
\ssA_{\pi,n}^{\Xi} \coloneqq \bigl\{ a\in \ssA_\pi \,\st \,\te_a \in \sspi_n \,, \he_a\in \Xi \bigr\}.
\]

Consider the characteristic vector $\one_n^\Xi \in\bV_\pi$ defined by
\[
\one_n^\Xi \coloneqq \one_{\ssA_{\pi,n}^\Xi} = \sum_{a\in \ssA_{\pi,n}^{\Xi}} \one_a. 
\]
Consider the vector subspace $\GlobSp_\pi \subseteq \RosSp_\pi$ defined by the \emph{global residue conditions}
\begin{equation}\label{eq:global}
\langle\one_n^\Xi,\cdot\rangle=0 \qquad \qquad \text{for all } \,\, 
n \in [r] \text{ and } \Xi \in \concom^*_{h< n}.
\end{equation}

Denote by $\Glob\subseteq\bV_\pi$ the collection of all the vectors $\one_n^\Xi$, that is,
\[
\Glob \coloneqq  \Bigl\{\one_n^\Xi \,\st\, n\in [r],\, \Xi\in \concom^*_{h< n}\Bigr\}.
\] 
This collection is obviously set-theoretically independent in $\bV_\pi$. We show the following result.

\begin{prop}\label{GRC = s-1} We have $\codim(\GlobSp_\pi,\RosSp_\pi)= \summit-\ssc$.
\end{prop}

For each connected component $\mC \in \concom_{h \leq n}$, denote by $\mC_n=\sspi_{n}\cap \rmC$ the set of vertices of level $n$ in $\mC$. Note that for each $n\in [r]$, 
each special component $\Xi\in \concom^*_{h < n}$ is contained in a unique component $\mC \in \concom_{h\leq n}$. 

We use a reasoning similar to the one that led to Equation~\eqref{eq:codim-Ros} to get
\begin{equation}\label{eq:codim-Glob}
\codim(\GlobSp_\pi,\bV^0_\pi)=
\sum_{n\in [r]}\, \sum_{\mC \in \concom_{h\leq n}}\codim\Bigl(\GlobSp_\pi\cap \bV_{\pi,\mC_{n}},\bV^0_{\pi,\mC_{n}}\Bigr).
\end{equation}
For each $n\in [r]$ and $\mC \in \concom_{h \leq n}$, define the three collections of characteristic vectors in $\bV_{\pi,\mC_{n}}$:
\begin{align*}
\LRC^{\,n,\rmC}&\coloneqq\bigl\{\,\one_v\,\st\, v\in \mC_n \textrm{ and }\one_v \neq 0\,\bigr\},\qquad \Ros^{\,n,\rmC}\coloneqq \bigl\{\,\one_e\,\st\, e \in\ssE_\pi^c(\mC_n)\,\bigr\}, \quad \text{ and }\\
\Glob^{\,n,\rmC}&\coloneqq \bigl\{\,\one_n^\Xi\,\st\, \Xi \in \concom_{h < n}^*\textrm{ and }\, \, \Xi\subseteq \mC \bigr\}.
\end{align*}
Clearly, each of $\Glob^{\,n, \rmC}\cup\Ros^{\,n, \rmC}$ and $\LRC^{\,n, \rmC}$ are set-theoretically independent collections of vectors in $\bV_{\pi,\mC_{n}}$. However, the two collections are related since
\begin{equation*}
   \bigcup_{\varphi \in \Glob^{n, \rmC}\cup\Ros^{n, \rmC}}[\varphi]= \E_{\rmC_{n}}\setminus \barA_\pi = 
   \bigcup_{\psi\in \LRC^{\, n, \rmC}}[\psi].
\end{equation*}

\begin{lemma}\label{lem:unrelated2}
Notation as above, for each $n \in [r]$ and $\mC \in \concom_{h\leq n}$, the two collections of vectors $\Glob^{\,n, \rmC}\cup\Ros^{\,n, \rmC}$ and $\LRC^{\,n, \rmC}$ in $\bV_{\pi,\rmC_{n}}$ are related but properly unrelated.
\end{lemma}

\begin{proof} 
We only need to show that the two collections are properly unrelated. Let $\Phi \subseteq \Glob^{\,n,\rmC}\cup\Ros^{\,n,\rmC}$ and $\Psi \subseteq \LRC^{\,n,\rmC}$ be two  nonempty subsets such that $[\Phi] = [\Psi]$, that is
\begin{equation}\label{TR}
   \bigcup_{\varphi \in \Phi}[\varphi] = \bigcup_{\psi\in \Psi}[\psi].
\end{equation}
If $\Psi=\LRC^{\,n,\rmC}$, then the right-hand side of \eqref{TR} is $\E_{\rmC_{n}}\setminus \barA_\pi$, which yields 
$\Phi=\Glob^{\,n,\Lambda}\cup\Ros^{\,n,\Lambda}$. 

Assume for the sake of contradiction that $\Psi \neq \LRC^{\,n, \mC}$. Let $Z$ be the proper subset of $\mC_n$ consisting of vertices $v$ in $\mC_{n}$ such that $\one_v\in\Psi$. 

If there were a horizontal edge $e$ connecting a vertex $v\in Z$ to a vertex $w\in \mC_n \setminus Z$, then the corresponding horizontal arrow $a=vw$ would be in $[\one_v]$, whence, in the right-hand side of \eqref{TR}. This would be possible only if $\one_{e}\in\Phi$, and so, $wv$ would belong to the left- and therefore, right-hand side of \eqref{TR}, yielding $w\in Z$, a contradiction.

Thus, we can assume there is no horizontal edge connecting $Z$ to $\mC_{n} \setminus Z$. Since $\mC$ is connected, there has to be an upward arrow $a = vu$ connecting a vertex $v$ in $Z$ to a vertex $u$ in a component $\Xi \in \concom_{h< n}^*$. Moreover, since $\mC$ is connected, $\Xi$ is contained in $\mC$, and there is a vertical arrow $b=wz$ connecting a vertex  $w \in \mC_{n} \setminus Z$ to a vertex $z$ in $\Xi$. Since $v\in Z$, the arrow $a$ appears in the right-hand side of \eqref{TR}, and thus, in the left-hand side as well. This means that $\one_{n}^\Xi \in \Phi$. 
We infer that  $b$ appears in the left-hand side of \eqref{TR}, and so in the right-hand side as well. This yields $w \in Z$, a contradiction.
\end{proof} 

\begin{lemma}\label{lem:card} For each $\mC \in \concom_{h \leq n}^*$, the following holds:
\begin{itemize}
\item $|\LRC^{\,n,\rmC}|=|\mC_{n}|$ unless $\mC$ is a singleton of level $n$, in which case $|\mC_{n}|=1$ but $|\LRC^{\,n,\rmC}|=0$.
\item $|\Ros^{\,n,\rmC}|$ is the number of horizontal edges in $\mC_{n}$.
\item $|\Glob^{\,n,\rmC}|$ is the number of connected components in $\concom^*_{h < n}$ that are included in $\rmC$. 
\end{itemize}
\end{lemma}

\begin{proof} The proof can be obtained by direct verification. We omit the details.
\end{proof}

We can now finish the proof of the proposition.

\begin{proof}[Proof of Proposition~\ref{GRC = s-1}] We will combine the previous lemmas with Proposition~\ref{prop:properlyunlinked}. Notice first that, for each $n\in [r]$,
\begin{equation}\label{eq:CC1}
\concom_{h< n} \setminus \concom^*_{h < n}=\concom_{h \leq n}\cap \concom_{h < n}.
\end{equation} 
It follows that $\mC_{n}=\emptyset$ if and only if $\mC\in \concom_{h < n}\setminus \concom^*_{h < n}$, in which case, 
$\Glob^{\,n, \rmC}$, $\Ros^{\,n, \rmC}$ and $\LRC^{\,n, \rmC}$ are all empty. In fact, these collections are all empty if and only if either, $\mC \in \concom_{h < n}\setminus \concom^*_{h < n}$ or, $\mC$ is a singleton whose only vertex has level $n$, in which case $\mC$ is an irreducible  summit of the level graph $G$. We apply Proposition~\ref{prop:properlyunlinked} and deduce from \eqref{eq:codim-Glob} and Lemma~\ref{lem:unrelated2} 
that 
\begin{equation}\label{eq:codG}
\codim\big(\GlobSp_\pi ,\bV_\pi \big)=\sum_{n\in [r]}\Big(\summit_{\ired}(n)+\sum_{\rmC \in \concom_{h \leq n}^*}\bigl(\,|\Glob^{\,n,\rmC}|+|\Ros^{\,n,\rmC}|+
|\LRC^{\,n,\rmC}|-1\big)\Big)
\end{equation}
where, for each $n\in [r]$, the quantity  $\summit_{\ired}(n)$ is the number of irreducible summits of the level graph $G$ whose unique vertex has level $n$, and
\begin{equation}\label{eq:CC2}
\concom_{h \leq n}^*=\concom_{h \leq n}\setminus \bigl(\concom_{h \leq n}\cap \concom_{h < n}\bigr).
\end{equation}

Now, applying Lemma~\ref{lem:card} to each $\mC \in \concom_{h \leq n}^*$,  we obtain the cardinalities in \eqref{eq:codG}, leading, for each $n \in [r]$, to
\begin{equation}\label{eq:Tele}
\sum_{\rmC\in \concom_{h \leq n}^*}\Big(|\Glob^{\,n,\rmC}|+|\Ros^{\,n,\rmC}|+
|\LRC^{\,n,\rmC}|-1\Big)=|\concom^*_{h < n}|+|E(\sspi_n)|+|\sspi_n|-\summit_{\ired}(n)-|\concom_{h \leq n}^*|.
\end{equation}
For each $n\in [r]$, by \eqref{eq:CC1} and \eqref{eq:CC2}, we have $|\concom^*_{h < n}|-|\concom_{h \leq n}^*|=|\concom_{h < n}|-|\concom_{h \leq n}|$. Therefore, adding the equations \eqref{eq:Tele} for all $n$, and using \eqref{eq:codG}, we deduce that
\[
\codim\big(\GlobSp_\pi,\bV_\pi\big)=-\ssc+|\ssE_\pi^c|+|V|.
\]
Finally, using Proposition~\ref{prop:LRC} and Proposition~\ref{codimRos}, we get
\[
\codim\big(\GlobSp_\pi,\RosSp_\pi\big)=-\ssc+|\ssE_\pi^c|+|V|-(|V|-\summit_{\ired})-(|\ssE_\pi^c|-\summit_{\red})=\summit-\ssc. \qedhere
\]
\end{proof}

\begin{remark}[The case of level graphs with a unique summit]
  Notice that $\codim\big(\GlobSp_\pi,\RosSp_\pi\big)=0$ if the level graph $(G,\pi)$
  has a unique summit (and is therefore connected). This is always the case if the underlying graph $G$ is (multi)complete, that is, each pair of vertices are connected by at least one edge. In this case, the global residue conditions are a consequence of the downward vanishing, local and Rosenlicht residue conditions. In the applications to the study of limits of spaces of Abelian differentials, we note that this is the setting of the work \cite{EM} and \cite{ES07}, in which global residue conditions do not play any role. 
\end{remark}

\subsection{Proof of Theorem~\ref{thm:dimG}} 
We are now in position to prove the dimension count. Proposition~\ref{prop:LRC} yields
\[\dim_{\k}\bV^0_\pi=2|\ssE_\pi^c|+|\ssE_\pi|-|V|+\summit_{\ired}.\]
Also, Proposition~\ref{codimRos} yields $\codim(\RosSp_\pi,\bV^0_\pi)=|\ssE_\pi^c|-\summit_{\red}$, whereas
Proposition~\ref{GRC = s-1} yields $\codim(\GlobSp_\pi,\RosSp_\pi)=\summit-\ssc$. Finally, $|E| = |\ssE_\pi| + |\ssE_\pi^c|$. Combining these all,  we get the desired formula $\dim_{\k}\GlobSp_\pi=|E|-|V|+\ssc=g(G)$, as required.
\qed

\subsection{An example} \label{sec:example-dimG} 
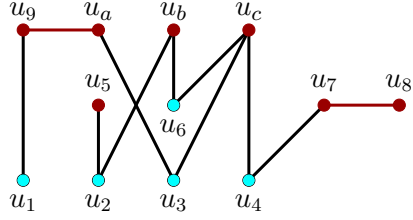
\begin{figure}
\centering
\begin{tikzpicture}
\draw[line width=0.4mm, darkred] (5,0) -- (6,0);
\draw[line width=0.4mm] (5,0) -- (5,-2);
\draw[line width=0.4mm] (6,0) -- (7,-2);
\draw[line width=0.4mm] (8,0) -- (7,-2);
\draw[line width=0.4mm] (6,-1) -- (6,-2);
\draw[line width=0.4mm] (7,0) -- (6,-2);
\draw[line width=0.4mm] (7,0) -- (7,-1);
\draw[line width=0.4mm] (8,0) -- (7,-1);
\draw[line width=0.4mm] (9,-1) -- (8,-2);
\draw[line width=0.4mm] (8,0) -- (8,-2);
\draw[line width=0.4mm, darkred] (10,-1) -- (9,-1);

\filldraw [darkred] (5,0) circle (0.8mm);
\filldraw [darkred] (6,0) circle (0.8mm);
\filldraw [darkred] (7,0) circle (0.8mm);
\filldraw [darkred] (8,0) circle (0.8mm);

\filldraw[darkred] (6,-1) circle (0.8mm);

\draw[line width=0.1mm] (7,-1) circle (0.8mm);
\filldraw[aqua] (7,-1) circle (0.7mm);

\filldraw[darkred] (9,-1) circle (0.8mm);

\filldraw[darkred] (10,-1) circle (0.8mm);

\draw[line width=0.1mm] (5,-2) circle (0.8mm);
\filldraw[aqua] (5,-2) circle (0.7mm);
\draw[line width=0.1mm] (6,-2) circle (0.8mm);
\filldraw[aqua] (6,-2) circle (0.7mm);
\draw[line width=0.1mm] (7,-2) circle (0.8mm);
\filldraw[aqua] (7,-2) circle (0.7mm);
\draw[line width=0.1mm] (8,-2) circle (0.8mm);
\filldraw[aqua] (8,-2) circle (0.7mm);

\draw[black] (5,0.25) node{$\ssu_9$};
\draw[black] (6,0.25) node{$u_a$};
\draw[black] (7,0.25) node{$u_b$};
\draw[black] (8,0.25) node{$u_c$};

\draw[black] (6,-.7) node{$\ssu_5$};

\draw[black] (7,-1.3) node{$\ssu_6$};

\draw[black] (9,-0.7) node{$\ssu_7$};
\draw[black] (10,-0.7) node{$\ssu_8$};

\draw[black] (5,-2.3) node{$\ssu_1$};
\draw[black] (6,-2.3) node{$\ssu_2$};
\draw[black] (7,-2.3) node{$\ssu_3$};
\draw[black] (8,-2.3) node{$\ssu_4$};

\end{tikzpicture} 
\caption{A level graph with three levels. Summits are drawn in red. There are three irreducible and two reducible summits. 
The subgraph $G[\{u_9,u_a\}]$ is a special connected component in $\concom_{h<3}$, but it is nonspecial in $\concom_{h<2}$.
}
\label{fig2}
\end{figure}

We discuss the example of the level graph $(G, \pi=(\sspi_1, \sspi_2, \sspi_3))$ with three levels depicted in Figure~\ref{fig2}. For simplicity, we let $ij$ denote the arrow from $u_i$ to $u_j$ for distinct $i,j$. Let $h\colon V \to [3]$ be the level function. 
For each $n\in [3]$, let 
\[
\LRC^{n}\coloneqq \bigl \{\one_v\,\st\,v\in\sspi_n\textrm{ and }\one_v\neq 0\bigr\},\quad 
\Ros^{n}\coloneqq \bigl \{\one_e\,|\,e\in E(\sspi_n)\bigr\},\quad
\Glob^{n}\coloneqq \bigl \{\one_n^\Xi\,\st \,\Xi\in \concom^*_{h< n}\bigr\}
\]
be subsets of $\bV_\pi$. Then, we obtain
\begin{align*}
\LRC^{\,3} &= \{\one_{19},\  \one_{25} + \one_{2b},\  \one_{3a} + \one_{3c},\ \one_{4c} + \one_{47}\},\qquad 
\Ros^{\,3}=\emptyset,\\ 
 \Glob^{\,3} &= \{\one_{19} + \one_{3a},\ \one_{25} ,\ \one_{2b} + \one_{3c} + \one_{4c},\ \one_{47}\},\\
 \LRC^{\,2} &= \{\one_{6b} + \one_{6c}, \ \one_{78},\ \one_{87}\},\qquad 
  \Ros^{\,2} = \{ \one_{78} + \one_{87}\},\qquad
\Glob^{\,2} = \{\one_{6b}, \ \one_{6c}\},\\
\LRC^{\,1} &= \{\one_{9a}, \ \one_{a9}\},\qquad
\Ros^{\,1} = \{\one_{9a} + \one_{a9}\},\qquad
\Glob^{\,1} = \emptyset.
\end{align*}
Notice that $G[\sspi_{3}]$ has four vertices (nonisolated in $G$) and no edges, and $G[\ssV_{h<3}]$ has four connected components. Hence, $\LRC^{\,3}$ and $\Glob^{\,3}$ have four elements each, whereas $\Ros^{\,3}=\emptyset$. The first two sets are related but properly unrelated, hence $\codim(\GlobSp_\pi\cap \bV_{\pi,\sspi_3}, \bV_{\pi,\sspi_3} ) = 7$; in other words, 
$\GlobSp_\pi\cap \bV_{\pi,\sspi_3}=0$. Also, $ \RosSp_\pi\cap \bV_{\pi,\sspi_3}=\bV_{\pi, \sspi_3}^0$, and  
$\codim(\RosSp_\pi\cap \bV_{\pi,\sspi_3},
\bV_{\pi,\sspi_3}) = 4$.

Further up, $G[\sspi_2]$ has four vertices, one isolated in $G[\ssV_{h\leq 2}]$, and one horizontal edge, whereas $G[\ssV_{h<2}]$ has three connected components, only two connected to a vertex in $\sspi_2$. Then, $\LRC^{\,2}$, $\Ros^{\,2}$ and $\Glob^{\,2}$ have respectively three, one and two elements. We have $\codim(\RosSp_\pi \cap \bV_{\pi,\sspi_2},
\bV_{\pi,\sspi_2}) = 3$, and $\codim(\GlobSp_\pi\cap \bV_{\pi,\sspi_2}, \bV_{\pi,\sspi_2} ) = 4$, the maximum possible.

Finally, $G[\sspi_1]$ has four vertices and one horizontal edge, whereas $G[\ssV_{h<1}]$ is empty. Then, $\LRC^{\,1}$, $\Ros^{\,1}$ and $\Glob^{\,1}$ have respectively two, one and zero elements. The codimension of 
$\GlobSp_\pi \cap \bV_{\pi,\sspi_1}=\RosSp_\pi\cap \bV_{\pi,\sspi_1}$ in $\bV_{\pi,\sspi_1}$ is 2, the maximum possible. 

It follows that $\dim_{\k}\GlobSp_\pi=0$ and $\dim_{\k}\RosSp_\pi=4$. Clearly, $\dim_{\k} \bV^0_\pi=4$. Notice that $G$ has $2$ horizontal arrows and five summits, two of which are reducible and three of them are irreducible, so $\summit_\red=2$, $\summit_\ired=3$, and $\summit=5$.

\section{Splitting, realization and degeneration} \label{sec:splitting}

Let $V$ be a finite nonempty set, and denote by $\ssub{2}!^V$ the family of subsets of $V$. We say a function $\eta \colon \ssub{2}!^V\to\R$ is \emph{nonnegative} if $\eta \geq 0$. We say $\eta$ is \emph{nonincreasing} if $\eta(J)\geq\eta(I)$ for $J\subseteq I\subseteq V$, and \emph{nondecreasing} if the reverse inequalities hold.

A function $\eta \colon \ssub{2}!^V\to\R$ is called \emph{submodular} if $\eta(\emptyset)=0$ and we have the inequalities
\[
  \eta(\ssI_1)+\eta(\ssI_2) \geq \eta(\ssI_1\cup \ssI_2)+ \eta(\ssI_1\cap \ssI_2)
\quad\text{for each }\ssI_1,\ssI_2\subseteq V.
\]
 The quantity $\eta(V)$ is called the \emph{range} of $\eta$. Similarly, $\chi \colon \ssub{2}!^V \to \R$ is called \emph{supermodular} if $\chi (\emptyset)=0$ and the inequalities above are all reversed; the quantity $\chi(V)$ is called the range of $\chi$. A function which is both submodular and supermodular is called \emph{modular}. Modular functions are in bijection with elements of $\R^V$: each $q\in \R^V$ can be viewed as a modular function $q\colon\ssub{2}!^V\to\R$ by setting  $q(I)\coloneqq \sum_{v\in I}q(v)$ for each $I\subseteq V$, with the convention that $q(\emptyset)\coloneqq 0$.

 For $\eta \colon \ssub{2}!^V\to\R$ with $\eta(\emptyset)=0$, we define the \emph{adjoint to $\eta$} denoted $\eta!^* \colon \ssub{2}!^V\to\R$ by
\[
\eta!^*(I)\coloneqq\eta(V)-\eta(V\setminus I)\quad\text{for each }I\subseteq V.
\]
It is easy to see that $\eta$ is submodular, resp.~supermodular, if and only if $\eta!^*$ is supermodular, resp.~submodular. Furthermore, $\eta$ and $\eta!^*$ have the same range, and $\ssub{(\eta!^*)}!^* = \eta$. If $\eta$ is submodular, then $\eta^*\leq\eta$. We refer to the pair $(\eta!^*, \eta)$ as a modular pair. Note that for a modular function $q$, we have $\ssub{q}!^*=q$, and so $(q,q)$ is a modular pair. 

To each submodular function $\eta\colon 2^V\to\R$ is associated a polytope,
\[
\Q_{\eta}\coloneqq\bigl\{q\in\R^V\,\st \, q(I)\leq\eta(I)\,\quad \forall I\subseteq V, \text{ with equality if }I=V\bigr\},
\]
the base polytope of the polymatroid associated to $\eta$, or simply, the (\emph{base}) \emph{polytope} of $\eta$. Equivalently, $\Q_\eta$ is defined as the set of points that verify the inequalities $q(I)\geq \eta!^*(I)$ for $I\subseteq V$, with equality for $I=V$.

\subsection{Submodular functions of subspaces} \label{sec:submodular-subspace} Let $\k$ be a field. For each $v\in V$, let $\VS_v$ be a vector space over $\k$. Let $\VS\coloneqq \bigoplus_{v\in V} \VS_v$. For each subset $I\subseteq V$, let $\VS_I\coloneqq \bigoplus_{v\in I}\VS_v$ and denote by $\proj{I}\colon \VS\to \VS_{I}$ the corresponding projection map.  

A finite-dimensional vector subspace $\ssW\subseteq \VS$ gives rise to a submodular function $\nu!^*_{\ssW}\colon \ssub{2}!^V \to \Z$ defined by setting $\nu!_{\ssW}^*(I) \coloneqq \dim_{\k}(\ssW_I)$ for each $I\subseteq V$ with $\ssW_I\coloneqq \proj{I}(\ssW)$. Denote by $\ssW^{I^c}$ the kernel of $\proj{I}\colon \ssW \to \ssW_I$, so that we have a short exact sequence 
\begin{equation}\label{eq:short-exact}
0 \to \ssW^{I^c} \to \ssW \to \ssW_{I} \to 0.
\end{equation}
Then, we have $\nu!_W(I) = \dim_{\k} \ssW^I$.

\begin{prop}\label{prop:submodular-subspace} The function $\nu!_{\ssW}^* \colon \ssub{2}!^V \to \Z$ is submodular, nonnegative and nondecreasing. 
\end{prop}

\begin{proof} 
From~\eqref{eq:short-exact}, we get $\ssW_{I} \simeq \rquot{\ssW}{\ssW^{I^c}}$, and so the adjoint $\nu!_W$ of $\nu!^*_{\ssW}$ is given by $\nu!_{\ssW}(I) = \dim_{\k}(\ssW^{I})$ for each $I\subseteq V$. It will be enough to show that $\nu!_{\ssW}$ is supermodular, that is, for $I, J \subset V$, we have 
\[
\dim_{\k} \ssW^{I}+\dim_{\k} \ssW^{J}\leq\dim_{\k} \ssW^{I\cup J}+\dim_{\k} \ssW^{I \cap J}.
\]
This follows from the observation that
$\ssW^{I}$ and $\ssW^{J}$ are subspaces of $\ssW^{I \cup J}$
with intersection $\ssW^{I \cap J}$. The rest is clear.
\end{proof}

\subsection{Filtration associated to an ordered partition} Let $\pi=(\sspi_1, \dots, \sspi_r)$ be an ordered partition of $V$, with the level function $h=\ssh_\pi \colon V \to [r]$.

Let $\filt_0 = \emptyset$, and for each $n \in [r]$, let $\filt_n \coloneqq \sspi_1 \cup \dots \cup \sspi_n =\ssV_{h\leq  n}$ be the set of all $v\in V$ with $h(v)\leq n$. We obtain an ascending filtration $\filt_{\bullet}$ indexed by $0, 1, \dots, r$, and the data of $\pi$ is equivalent to that of $\filt_\bullet$.

\subsection{Splitting}
Let $\pi=(\sspi_1, \dots, \sspi_r)$ be an ordered partition of $V$, and $\filt_\bullet$ and $\rmH_\bullet$ the corresponding ascending filtrations associated to $\pi$, respectively. 

To each supermodular function $\mu\colon \ssub{2}!^V\to\R$ we associate another supermodular function $\mu!_{\pi}$ called the \emph{splitting} of $\mu$ with respect to $\pi$, as follows.

First, for each $n\in[r]$, we define the function  
\[
\mu!_{\filt_n/\filt_{n-1}}\colon \ssub{2}!^V\to\Z
\]
by setting 
\[
\mu!_{\filt_n/\filt_{n-1}}(I)\coloneqq \mu((I\cap\filt_n)\cup\filt_{n-1})-\mu(\filt_{n-1})\quad\text{for each } I\subseteq V.
\]
This is a supermodular function, see \cite[\S2.4]{AE-modular-polytopes}. 

Then, we define $\mu!_\pi$ by
\[
\mu!_{\pi}\coloneqq\sum_{n=1}^r\mu!_{\filt_n/\filt_{n-1}}.
\]

\subsection{Realization} We now show that if the supermodular function $\mu$ is associated to a subspace, then any of its splittings is also associated to a subspace. More precisely, let $\pi=(\sspi_1, \dots, \sspi_r)$ be an ordered partition of $V$ and denote by $\filt_\bullet$ the corresponding filtration of $V$. Notation as in \ref{sec:submodular-subspace}, consider a subspace $\ssW\subseteq\VS$. We associate to $\pi$ a subspace $\ssW(\pi)\subseteq\VS$ defined as follows.

Let
\[
0 \subseteq \ssW^{\filt_1} \subseteq \ssW^{\filt_2} \subseteq \dots \subseteq \ssW^{\filt_n}=\ssW.
\]

For each $n\in[r]$, define $\ssW_n$ as the quotient $\rquot{\ssW^{\filt_n}}{\ssW^{\filt_{n-1}}}$ so that we have the following short exact sequence 
\begin{equation}\label{eq:short-exact2}
0 \to \ssW^{\filt_{n-1}} \to \ssW^{\filt_n} \to \ssW_{n} \to 0.
\end{equation}
Note that, we have $\ssW_n = \bigl(\ssW^{\filt_n}\bigr)_{\sspi_n} = \proj{\sspi_n}(\ssW^{\filt_n})$. We define
\[
\ssW(\pi)\coloneqq\bigoplus_{n=1}^r \ssW_n.
\]
Since each $\ssW_n$ lives naturally in $\VS_{\sspi_n}=\bigoplus_{v\in \sspi_n} \VS_v$, by using $\VS \simeq \bigoplus_{n\in[r]}\VS_{\sspi_n}$, we view  
\[
\ssW(\pi) \subseteq \VS =\bigoplus_{v\in V} \VS_v.
\]
\begin{prop}\label{prop:supermodular2} Let $\nu=\nu!_\ssW$ be the supermodular function of a subspace $\ssW \subseteq \VS$. For each ordered partition $\pi$ of $V$,  the splitting  $\nu_{\pi}$ is the supermodular function associated to $\ssW(\pi) \subseteq \VS$.
\end{prop}

\begin{proof} For each $I\subseteq V$, we have the following commutative diagram 
\[
\begin{CD}
   0 @>>> \ssW^{\filt_{n-1}} @>>> \ssW^{\filt_n} @>>> (\ssW^{\filt_n})_{\sspi_n}@>>>0\\
   @. @VVV@VVV @VVV @.\\
   0 @>>> \ssW^{\filt_{n-1}\cap(I\cup\sspi^c)} @>>> \ssW^{\filt_n} @>>> (\ssW^{\filt_n})_{\ssI^c \cap \sspi_n} @>>> 0
\end{CD}
\]
with each of the two rows forming a short exact sequence.  Using the above diagram, and the identity $\ssW_n = \bigl(\ssW^{\filt_n}\bigr)_{\sspi_n}$, we deduce the equality
\[
\nu!_{\ssW_n}(I) = \nu(\filt_n)-\nu(\filt_{n-1}) - \Bigl(\nu(\filt_n) -  \nu\bigl(\filt_{n-1} \cup (\filt_n\cap I)\bigr) \Bigr) = \nu!_{\filt_n/\filt_{n-1}}(I) \qquad \text{for each } I\subseteq V.
\]
Therefore, we have $\nu!_{\ssW_n} = \nu!_{\filt_n/\filt_{n-1}}$, and 
\[
\nu!_{\ssW(\pi)} = \sum_{n=1}^r \nu!_{\ssW_n} = \sum_{n=1}^r \nu!_{\filt_n/\filt_{n-1}} = \nu!_\pi,
\]
as required.
\end{proof}

\subsection{Degeneration}\label{sec:upper-semicontinuity} 
An alternate point of view on $\ssW(\pi)$ is through degenerations, as follows.

First, denote by $\Gm$ the multiplicative group with $\Gm(\k) = \k^\times = \k\setminus\{0\}$.
The group $\Gm^V(\k) = \ssub{\bigl(\k^\times\bigr)}!^V$ acts componentwise on $\VS$.

Let $\varR = \k[t,t^{-1}]$ be the ring of Laurent polynomials with coefficients in $\k$. Let $\VS_\varR= \VS\otimes_{\k} \varR$. The group $\Gm^V(\varR)$ acts componentwise on $\VS_\varR$.

For each integer-valued level function $d\colon V\to\Z$, $v\mapsto \ssd_v$, with induced ordered partition given by $\pi$, let $x\in\Gm^V(\varR)$ be given by $\ssx_v\coloneqq t^{-\ssd_v}$ for each $v\in V$. Consider the multiplication by $x$ in $\VS_\varR$, $\varphi \mapsto x\cdot \varphi$, given by
\[
\ssub{(x\cdot\ssvarphi)}_{v}\coloneqq
\ssx_{v}\ssvarphi_v=t^{-\ssd_v}\ssvarphi_v \qquad \text{for each $\varphi\in\VS_\varR$ and each $v\in V$}.
\]
 Given a subspace $\ssW \subseteq \VS$, we obtain a submodule$x\cdot \ssW_\varR \subseteq \VS_{\varR}$ for $\ssW_\varR\coloneqq \ssW\otimes_{\k}\varR$. Then, evaluation at $t=\lambda \in \k\setminus\{0\}$ gives a subspace $x(\lambda)\cdot \ssW \subset \VS$. We denote by $\lim_{t\to 0}x\cdot\ssW$ the limit of the subspaces $x(\lambda)\cdot \ssW \subset \VS$ for $\lambda \neq 0$ in the Grassmannian $\grass(m, \VS)$ of $m$-dimensional subspaces of $\VS$ with $m = \dim_{\k}\ssW$.

\begin{prop}\label{prop:splitting} We have
\[
\lim_{t\to 0}x\cdot\ssW=\ssW(\pi).
\]
\end{prop}

\begin{proof} To simplify, let $\ssW^n = \ssW^{\filt_n}$. For each $i=1,\dots,r$, let $\ssd_n$ be the value taken by $d$ on $\sspi_n$. We have $\ssd_1<\ssd_2<\cdots\ssd_r$. 

Let $\ssm_n\coloneqq \dim\ssW_n$ for $n\in[r]$. Recall that $\ssW_n = \bigl(\ssW^{\filt_n}\bigr)_{\sspi_n} \subseteq \VS_{\sspi_n}$.

For each $n\in[r]$, we choose a lifting $\ssy_{n,1},\dots,\ssy_{n,\ssm_i}$  to $\ssW^{n}$ of a basis  $\ssz_{n,1},\dots,\ssz_{n,\ssm_{i}}$ of $\ssW_n $. Let $\ssz_{n,1},\dots,\ssz_{n,\ssm_n},\dots,\ssz_{n,\ssM_n}$ be an extension of the basis of $\ssW_n$ to one of $\VS_{\sspi_n}$. Order the basis $\ssz_{n,j}$ of $\VS$ lexographically, so that $\ssz_{i,j} < \ssz_{i,j}$ if either $i<i'$, or $i=i'$ and $j<j'$. Extend this order to one for the basis obtained from the $\ssz_{i,j}$ for the $m$-th exterior product $\bigwedge^m\VS$, where $m = \ssm_1+\dots+\ssm_r=\dim_{\k}\ssW$. 

The $\ssy_{n,j}$ for $j\leq \ssm_n$ form a basis of $\ssW$. Their exterior product $\ssy_{1,1}\wedge \ssy_{1,2}\wedge\cdots\wedge\ssy_{r,\ssm_{r}}$ can be written as the exterior product $\ssz_{1,1}\wedge\ssz_{1,2}\wedge\cdots\wedge \ssz_{r,\ssm_{r}}$ plus a linear combination of terms of higher order (with respect to the lexicographic order). 

Likewise, the $x\cdot\ssy_{i,j}$ form a basis of $x\cdot\ssW_\varR$ over the ring $\varR$, and we have
\[
(x\cdot\ssy_{1,1})\wedge\cdots\wedge(x\cdot\ssy_{r,\ssm_{r}})=
t^{-\sum_n \ssm_n\ssd_n}\ssz_{1,1}\wedge\ssz_{1,2}\wedge\cdots\wedge \ssz_{r,\ssm_{r}}+t^l z'
\]
for some $z'$ in $\bigwedge^m\VS_\varR$, where $l>-\sum_n \ssm_n\ssd_n$. Clearly, the limit of $(x\cdot\ssy_{1,1})\wedge\cdots\wedge(x\cdot\ssy_{r,\ssm_{r}})$ as $t$ approaches $0$ is  $\ssz_{1,1}\wedge\ssz_{1,2}\wedge\cdots\wedge \ssz_{r,\ssm_{r}}$. Since the $\ssz_{n,j}$ for $j\leq \ssm_n$ form a basis of the sum $\ssW(\pi)=\bigoplus_n\ssW_n$, the proof is complete.
\end{proof}

\section{Proof of Theorem~\ref{thm:faces-residue-polytope}}\label{sec:respoly} 

Consider the trivial level structure $\sspi_0$ on $G$ given by the ordered partition of $V$ into a unique set $V$. In this case, vanishing along downward arrows and global residue conditions are vacuum. The residue space $\GlobSp_{\sspi_0} \subset \bV = \bigoplus_{v\in V}\bV_v$ is the first homology of $G$ with $\k$-coefficients, and the polytope associated to $\nu!^*_{\GlobSp_{\sspi_0}}$ is the residue polytope of $G$  denoted by $\RP=\RP(G)$.

Since $\dim\GlobSp_{\sspi_0} = g(G)$ and both $\nu!_{\GlobSp_{\sspi_0}}$ and $\nu!^*_{\GlobSp_{\sspi_0}}$ are nonnegative, the residue polytope $\RP$ lives in the standard simplex $\Delta!_{g} \subset \R_{\geq 0}^V$, consisting of points whose coordinates sum up to $g$.

For the ordered partition $\pi = (\sspi_1, \dots, \sspi_r)$, consider the subspace $\GlobSp_{\pi}\subseteq\Upsilon$. Denote by $\gamma!_{\pi}=\nu!^*_{\GlobSp_{\sspi}}$ the  corresponding submodular function and by $\P_{\pi}$ the corresponding polytope, called the $\pi$-residue polytope. Again, we have  $\P_{\pi} \subseteq \Delta!_{g}$.

Consider a coarsening $\pi'=(\sspi_1', \dots, \sspi_s')$ of $\pi$, that is, an ordered partition of $V$ such that the natural surjection $h'\colon \ssV \to [s]$ factors through $h\colon \ssV\to [r]$, and the induced map  $c\colon[r]\to [s]$ satisfies $c(\ssn_1)\leq c(\ssn_2)$ for each $\ssn_1,\ssn_2\in[r]$ with $\ssn_1\leq \ssn_2$. 
Notation as in the previous section, the group $\Gm^V(\varR)$ acts componentwise on $\bV_\varR = \bigl(\bigoplus_v\bV_v\bigr)\otimes_{\k}\varR$, with $\bV_v = \k^{\E_v}$.

Consider a function $d\colon V\to\Z$, $v\mapsto \ssd_v$, with induced ordered partition given by $\pi$, and let $x\in\Gm^V(\varR)$ be given by $\ssx_v\coloneqq t^{-\ssd_v}$ for each $v\in V$. Mutiplication by $x$ is given by 
\[
\ssub{(x\cdot\ssvarphi)}_{a}\coloneqq
\ssx_{v}\ssvarphi_a=t^{-\ssd_v}\ssvarphi_a \qquad \text{for each $\varphi\in\bV$, each $v\in V$ and each $a\in\E_v$}.
\]
Applying Proposition~\ref{prop:splitting} to the subspace $\GlobSp_{\pi'} \subset \bV$, we get
\begin{equation}\label{eq:splitting2}
\lim_{t\to 0}x\cdot\GlobSp_{\pi'}=\GlobSp_{\pi'}(\pi).
\end{equation}

\begin{thm}\label{thm:limGlobSp} We have 
\[
\lim_{t\to 0}x\cdot\GlobSp_{\pi'}=\GlobSp_{\pi}.
\]
In particular, we have $\GlobSp_{\pi}=\GlobSp_{\pi'}(\pi)$, and so  $\gamma!_{\pi}$ is the splitting of $\gamma!_{\pi'}$ with respect to $\pi$. As a consequence, we have $\gamma_\pi \leq \gamma_{\pi'}$.
\end{thm}

\begin{proof} Dividing by the appropriate power of $t$ and then putting $t=0$ each of certain equations describing $x\cdot\RosSp_{\pi'}$, we get already some of the equations that describe the limit of $x\cdot\GlobSp_{\pi'}$, to wit, vanishing along downward arrows in $\barA_\pi$, local residue conditions, and Rosenlicht conditions for $\pi$:
\begin{align*}
\ssvarphi_{a}=0&\qquad\qquad \text{for each downward }a\in\bar\ssA_\pi,\\
\sum_{a\in\E_v}\ssvarphi_a=0&\qquad\qquad \text{for each }v\in V,\\
\ssvarphi_a+\ssvarphi_{\bar a}=0&\qquad\qquad \text{for each } a\in\E \setminus \E_\pi \text{ (horizontal for }\pi).
\end{align*}

Consider $\sspi'_i$  for $i\in [s]$. For each $i\in [s]$, let $\lvl_{i,1},\dots,\lvl_{i,q_i}$ be the parts of $\pi$ which are contained in $\sspi'_i$ in ascending order, and $\ssd_{i,1},\dots,\ssd_{i,q_i}$ the ascending sequence of levels attributed by $d$ to each of them. The order given by $\pi$ on the $\lvl_{i,j}$, $i\in [s]$, $j\in[\ssq_i]$, is thus the lexographic on the subindices. Let $\ssV_{h' < i}$ denote the union of the $\sspi_{n}'$ for $n< i$ in $[s]$, and let $\ssV_{h \subface_{\rm{lex}} (i,j)}$ denote the union of the $\lvl_{n,m}$ for $(n,m)\subface_{\rm{lex}}(i,j)$. We need to show that, for each $i\in [s]$ and $j=1,\dots,\ssq_i$, the limit of $x\cdot\GlobSp_{\pi'}$ satisfies the global residue conditions, namely, the equations
\begin{equation}\label{eq:missingequations1}
\sum_{a\in\E(\lvl_{i,r}, \Xi)} \ssvarphi_a=0\qquad\qquad \textrm{for each connected component }\Xi\textrm{ of }G[\ssV_{\subface (i,r)}].
\end{equation}

Now, for the component $\Xi$ of $G[\ssV_{h \subface_{\rm{lex}} (i,j)}]$, let $S_l$ be the subset of vertices of $\lvl_{i,l}$ it contains for each $t=1,\dots,j-1$. Then, $\Xi$ is the union of certain connected components $\Xi_1,\dots,\Xi_k$ of $G[\ssV_{h'< i}]$ to which we add the vertices in the $S_l$ for $l=1, \dots, j-1$, and the edges connecting the vertices in these sets. The arrows that connect $\lvl_{i,j}$ to $\Xi$ are all upward and have head either in $S_l$ for some $l=1, \dots,j-1$ or in $\Xi_l$ for some $l=1,\dots, k$. We get Equation \eqref{eq:missingequations1} summing up the following series of equations satisfied by $x\cdot\GlobSp_{pi'}$ for $t\neq 0$, by multiplying by $t^{-\ssd_{i,j}}$ and then putting $t=0$:
\begin{align*}
\sum_{l=1}^k\sum_{j=1}^{\ssq_i} t^{\ssd_{i,j}} \sum_{a\in\E(\lvl_{i,s},\Xi_l)} \ssvarphi_a=&0 &\textrm{(Global residue conditions for $\pi'$)}\\
-\sum_{l=1}^{j-1} t^{\ssd_{i,l}} \sum_{a\in\E_{S_l}}\varphi_a=&0
&\textrm{(Local residue conditions)}\\
\sum_{l=1}^{j-1} \sum_{j'=1}^{\ssq_i} \sum_{a\in\E(\ssS_l, \lvl_{i,j'})}  \big(t^{\ssd_{i,l}}\ssvarphi_a+t^{\ssd_{i,j'}}\ssvarphi_{\bar a}\big)=&0 
&\textrm{(Rosenlicht conditions)}\\
\sum_{l=1}^{j-1}\sum_{a\in\E(\ssS_l, \ssV_{h'> i})}t^{\ssd_{i,l}}\ssvarphi_a=&0. 
&\textrm{(Vanishing along downward arrows for $\pi'$)}
\end{align*}
Having established that the limit of $x\cdot\GlobSp_{\pi'}$ as $t$ approaches $0$ satisfies all the equations that $\GlobSp_\pi$ does finishes the proof of the first statement, as the two spaces $x\cdot\GlobSp_{\pi'}$ and $\GlobSp_\pi$ have the same dimension. 

Combining this with \eqref{eq:splitting2}, we get $\GlobSp_{\pi}=\GlobSp_{\pi'}(\pi)$, which implies that $\gamma!_{\pi}$ is the splitting of $\gamma!_{\pi'}$ with respect to $\pi$.

The last assertion follows from the splitting statement, or directly from the limit, using the lower semicontinuity of the dimensions of projections in the limit.
\end{proof}

 We can now prove Theorem~\ref{thm:faces-residue-polytope}.  
\begin{proof}[Proof of Theorem~\ref{thm:faces-residue-polytope}] We show that for each ordered partition $\pi$, the polytope $\P_{\pi}$ is a face of the residue polytope $\RP$. Moreover, each face of $\RP$ is of this form.

By \cite[Prop.~2.7]{AE-modular-polytopes}, the faces of $\RP(G)$ are the base polytopes of the splittings $\mu!_{\pi}$ of the supermodular function $\mu\coloneqq\nu!_{\GlobSp_{\pi_0}}$ with respect to the ordered partitions $\pi$ of $V$; see \cite[\S2.4]{AE-modular-polytopes}. But, $\mu!_{\pi}$ is the supermodular function of $\GlobSp_{\pi_0}(\pi)$ by \eqref{eq:splitting2}, and by Theorem~\ref{thm:limGlobSp}, we have $\GlobSp_{\pi_0}(\pi) = \GlobSp_{\pi}$.  We conclude that the base polytope of $\mu!_\pi$ coincides with that of $\GlobSp_\pi$, and the theorem follows. 
\end{proof}

\bibliographystyle{alpha}
\bibliography{bibliography}
\end{document}